\documentclass{stile}

\address{ }

\email{delpadro@dima.unige.it}
\urladdr{www.dima.unige.it/\~{}delpadro}

\usepackage{amsmath}
\usepackage{amsthm}
\usepackage{amssymb} 
\usepackage{lscape}

\newcommand{\BibTeX}{{\scshape Bib}\kern-.08em\TeX}
\newcommand{\T}{\S\kern .15em\relax}

\theoremstyle{plain}

\newtheoremstyle{newremark}
  {2.0pt}
  {2.0pt}
  {\rmfamily}
  {}
  {\rmfamily\bfseries}
  {.}
  {.5em}
  {}

\theoremstyle{definition}

\theoremstyle{remark}

\DeclareMathAlphabet{\mathbbb}{U}{bbold}{m}{n}
\DeclareMathAlphabet{\mathmanual}{U}{manfnt}{m}{n}

\newcommand{\calA}{\mathcal{A}}

\newcommand{\Sim}{\Sigma}

\newcommand{\Hom}{\mathrm{Hom}}

\newcommand{\Q}{\mathbb{Q}}

\newcommand{\tens}{\otimes}

\def\c#1{\mathop{ {\mathcal #1}  }\nolimits}

\def\ch{\mathop{\mathcal Ch}\nolimits}
\def\cch{\mathop{\mathcal CCh}\nolimits}

\def\ch-{\mathop{{\mathcal Ch}_-}\nolimits}
\def\cch-{\mathop{{\mathcal CCh}_-}\nolimits}
\def\ch+{\mathop{{\mathcal Ch}_+}\nolimits}
\def\cch+{\mathop{{\mathcal CCh}_+}\nolimits}

\def\cid#1{\mathop{{\rm Id}_{#1}}\nolimits}

\def\dual{{\displaystyle\check{\;}}}

\def\NN{{\mathbb N}}
\def\ZZ{{\mathbb Z}}
\def\QQ{{\mathbb Q}}

\def\CC{{\mathbb C}}
\def\FF{{\mathbb F}}

\def\KK{{\mathbb K}}

\def\fraz#1 #2 {\frac{#1}{#2}}

\def\h#1{{\kern .07em}^h\kern-.04em{#1}}
\def\a#1{{\kern .07em}{#1}_a}

\let\sem=\bf

\def\To{\longrightarrow}

\newcommand{\BD}{\begin{definition}}
\newcommand{\BDS}{\begin{definitions}}
\newcommand{\ED}{\end{definition}}
\newcommand{\EDS}{\end{definitions}}

\newcommand{\BR}{\begin{remark}}
\newcommand{\BRS}{\begin{remarks}}
\newcommand{\ER}{\end{remark}}
\newcommand{\ERS}{\end{remarks}}

\newcommand{\BN}{\begin{notation}}
\newcommand{\BNS}{\begin{notations}}
\newcommand{\EN}{\end{notation}}
\newcommand{\ENS}{\end{notations}}

\newcommand{\BE}{\begin{example}}
\newcommand{\BES}{\begin{examples}}
\newcommand{\EE}{\end{example}}
\newcommand{\EES}{\end{examples}}

\newcommand{\BT}{\begin{theorem}}
\newcommand{\ET}{\end{theorem}}

\newcommand{\BP}{\begin{proposition}}
\newcommand{\EP}{\end{proposition}}

\newcommand{\BC}{\begin{corollary}}
\newcommand{\EC}{\end{corollary}}

\newcommand{\BL}{\begin{lemma}}
\newcommand{\EL}{\end{lemma}}

\newcommand{\Bel}{\begin{description}}
\newcommand{\Eel}{\end{description}}

\def\c#1{\mathop{ {\mathcal #1}  }\nolimits}

\def\ch{\mathop{\mathcal Ch}\nolimits}
\def\cch{\mathop{\mathcal CCh}\nolimits}

\def\ch-{\mathop{{\mathcal Ch}_-}\nolimits}
\def\cch-{\mathop{{\mathcal CCh}_-}\nolimits}
\def\ch+{\mathop{{\mathcal Ch}_+}\nolimits}
\def\cch+{\mathop{{\mathcal CCh}_+}\nolimits}

\def\cid#1{\mathop{{\rm Id}_{#1}}\nolimits}

\def\dual{{\displaystyle\check{\;}}}

\def\NN{{\mathbb N}}
\def\ZZ{{\mathbb Z}}
\def\QQ{{\mathbb Q}}

\def\CC{{\mathbb C}}
\def\FF{{\mathbb F}}
\def\KK{{\mathbb K}}
\def\PP{{\mathbb P}}

\DeclareMathAlphabet{\mathbbb}{U}{bbold}{m}{n}
\newcommand{\I}{\mathbbb{1}}

\def\fraz#1 #2 {\frac{#1}{#2}}

\def\h#1{{\kern .07em}^h\kern-.04em{#1}}
\def\a#1{{\kern .07em}{#1}_a}

\def\Tr{{\mathsf{Tr} }}

\def\cp{{\mathsf{cp} }}

\let\sem=\bf

\def\To{\longrightarrow}

\let\iff=\Leftrightarrow

\usepackage[vcentermath]{youngtab}

\title{Integral objects and Deligne's category $\textrm{Rep}(S_{t})$.}
\author{Alessio Del Padrone}

\begin{document}

\maketitle

\begin{abstract}
We give negative answers to certain questions on abelian semisimple $\otimes$-categories
raised by Kahn and Weibel in connection with the preprint of Kahn
``On the multiplicities of a motive''. For the most interesting examples we used Deligne's 
category $\textrm{Rep}(S_{t},\FF)$ of representations of  the ``symmetric group $S_{t}$ with $t$ not an integer''
with $\FF$ any algebraically closed field of characteristic zero.
\end{abstract}

\section*{Introduction}
In the preprint \cite{KaMM} (now published as Parts I and IV of 
\cite{KaZFM}) Bruno Kahn, studying the rationality
of certain zeta functions,
semisimple rigid $k$-linear $\otimes$-category $\c{A}$ and
defined the notion of {\it object (geometrically) of integral type}. 
This is quickly reviewed, together with some basic terminology, in sections 1 and 2. 
His work raised some questions of general interest for  $\otimes$-categories.
In this note we give examples showing that:
\begin{enumerate}
\item there are abelian semisimple geometrically integral categories with non Schur-finite objects (see sections 3 and \ref{NotSF})
 answering a question of Weibel to Kahn (\cite[Remark 2.3]{KaMM} and \cite[Remark 2.3]{KaZFM}),

\item the $\otimes$-product (even a tensor square) of geometrically integral objects in an abelian semisimple  category need not be of integral type (see section 7), answering a question of Kahn to the author.
\end{enumerate}

In section 3 we give a ``toy example'' based on the free rigid tensor category on one object. Then, in section 4, we review
the main tool which is Deligne's category $\textrm{Rep}(S_{t},\FF)$ of representations of  the ``symmetric group $S_{t}$ with $t$ not an integer'', with $\FF$ any algebraically closed field of characteristic zero.
This is a $\otimes$-category, ``new'' in some sense (\cite[Introduction]{DeSt04}),  with interesting properties and its study is of independent interest. 
In Proposition \ref{repstnotschur} we give two proofs of the fact that the canonical generator of $\textrm{Rep}(S_{t},\FF)$ is not Schur-finite for $t\in\FF\setminus\NN$ (that is, for such a $t$, the length of the tensor powers of its canonical $\otimes$-generator grows more than exponentially, see \cite[Proposition 0.5 (i)]{DeCT02}). The first proof uses the universal property of $\textrm{Rep}(S_{t},\FF)$, studied in section 5, while the second proof is more direct.
This gives also an algebraic proof that for any such $t$ the category $\textrm{Rep}(S_{t},\FF)$ is not $\otimes$-equivalent to a category of super-representation of a super group scheme (cf. \cite[Th\'eor\`eme 0.6]{DeCT02} and \cite[Introduction]{DeSt04}), because there are not even $\otimes$-functors from it to the $\otimes$-category of supervector spaces (see Corollary \ref{NoSuper}).

We also incidentally note that Deligne's category $\c{A}=\textrm{Rep}(S_{t},\FF)$ disproves the claims \cite[4.4, 4.5]{HoHai}, which was not previously noticed.
As Professor Deligne kindly remarked in a letter to the author (22 mars 2007) also those in his \cite{DeSEGL}, as well as \cite[(1.27)]{DMOS}, disprove the claims.

\section{Terminology}

\subsection{$\otimes$-categories and $\otimes$-functors}
In what follows $\c{A}$ denotes a {\it $\otimes$-category}, by this I mean that $\c{A}$ is an additive pseudoabelian, {\it i.e.} Cauchy complete (see \cite[6.5]{Bo1} and \cite[I.1]{AK}), (strict) monoidal category (as defined in \cite[VII, p.162]{MacLCWM}, see also \cite[XI.3, Th. 1]{MacLCWM}) where the monoidal structure $-\otimes-\colon \c{A}\times\c{A}\To\c{A}$ is a biadditive functor. We denote by $\I$ its tensor unit, so the (necessarily) commutative endomorphism monoid $\c{A}(\I)$ is actually a commutative unitary ring. 
If $R$ is a commutative unitary ring, we say that $\c{A}$ is $R$-linear if $\c{A}(\I)$ is an $R$-algebra, in this case all $\c{A}(X,Y)$ are $R$-modules and $-\otimes-\colon \c{A}\times\c{A}\To\c{A}$ is in fact $R$-bilinear. Note that $\c{A}$ is always $\c{A}(\I)$-linear. 

A {\it $\otimes$-functor} (called also {\it strong monoidal functor} by Mac Lane) $F\colon\c{A}\To\c{B}$ between two such $\otimes$-categories is a functor equipped with families of natural transformations expressing compatibility conditions such as $F(\I_{\c{A}})\cong \I_{\c{B}}$ and $F(X\otimes Y)\cong F(X)\otimes F(Y)$ for any $X,Y$. We refer to \cite[XI.2]{MacLCWM},
\cite{DMOS}, or \cite{SaaCT} for precise definitions and the neeeded commutative diagrams.

\subsection{Symmetry}
We assume, without explicit mention, that each $\otimes$-category $\c{A}$ used in this paper is {\it symmetric}, that is we are also given a family of natural isomorphisms
$
\tau_{X,Y}\colon X\otimes Y\To Y\otimes X$ such that $\tau_{X,Y}^{-1}=\tau_{Y,X}$ for each $X,Y\in\c{A}$. 
Each {\it $\otimes$-functor} is therefore also required to respect the symmetry (cf. \cite[XI.2(10)]{MacLCWM}).

\subsection{Rigidity}
We say that such an $\c{A}$ is {\it rigid} if for each object $X$ there are morphisms
$
\varepsilon_{X}\colon X\dual\otimes X\To \I\mbox{ (evaluation), }
\eta_{X}\colon \I\To X\otimes X\dual\mbox{ (coevaluation) }
$
satisfying the following ``triangular identities''
$$
\cid{X}=(X=\I\otimes X\xrightarrow{\eta_{X}\otimes X} X\otimes X\dual \otimes X\xrightarrow{X\otimes \varepsilon_{X}}X\otimes \I=X)
$$
and 
$$
\cid{X\dual}=(X\dual=X\dual\otimes\I\xrightarrow{X\dual\otimes \eta_{X}} X\dual\otimes X \otimes X\dual\xrightarrow{\varepsilon_X\otimes X\dual}\I\otimes X\dual=X\dual).
$$
The object $X\dual$ is then  called the {\it (rigid) dual} of $X$ (see for example \cite[II.6]{AK} or \cite[pag. 72]{JoSt} for more details).

\subsection{Trace and Euler characteristic}
Rigid $\otimes$-categories are (canonically) {\it traced} (as defined in \cite{JSV}) by
$
\Tr=\Tr^{X}_{\I,\I}\colon\c{A}(X)\To\c{A}(\I), \Tr(f):=
\varepsilon_{X}\circ\tau_{X,X\dual}\circ (f\otimes X\dual)\circ\eta_{X}.
$
The {\it Euler characteristic}, $\chi(X)$, of an object $X$ (a.k.a. its {\it rigid dimension}) is then the categorical trace of its identity, that is $\chi(X):=\Tr(\cid{X})\in \c{A}(\I)$. 

More generally the trace $\Tr^{U}_{A,B}(f)\colon A\To B$
of a morphism $f\colon A\otimes U\To B\otimes U$ is defined as
$$
\Tr^{U}_{A,B}(f)=(B\otimes(\varepsilon_{U}\circ\tau_{U,U\dual}))\circ (f\otimes U\dual)\circ (A\otimes \eta_{U}).
$$ 

\subsubsection{{\sem $\otimes$-functor and traces}} Rigidity is obviously preseved by any $\otimes$-functor $F\colon\c{A}\To\c{B}$, 
hence the canonical traces are also preserved by such a functor. In particular $\chi(F(X))=F(\chi(X))$. If moreover
$F$ is $\c{A}(\I)$-linear then $\chi(X)\cid{I_{\c{B}}}=F(\chi(X))=\chi(F(X))$.

The free traced monoidal category has been described in \cite{JSV} and \cite{Abr}. 
Thanks, for example, to \cite[Proposition 3]{Abr} we have the following proposition.

\subsubsection{{\sem Proposition}}\label{Abramsky}
{\it
Let $\c{A}$ be a traced monoidal category with zero object,
$A$ an object of $\c{A}$, $n\in \NN_+$ and $\sigma\in \Sigma_n$.
\begin{enumerate}
\item
$A^{\otimes n}=0$ if and only if $A=0$.
\item
$
\Tr(\sigma_{A,\dots,A}\colon A^{\otimes n}\to A^{\otimes n})=
{\chi(A)}^{|\mbox{\rm{cycles of } }\sigma|}.
$
\end{enumerate}
}

\subsection{The monoidal ideal $\c{N}$}\label{NumTrivial}
In a traced category with zero object $\c{A}$ we have the sets of 
{\it morphisms universally of trace zero}
$$
\c{N}(X,Y):=\{g\in\c{A}(X,Y)\mid  \Tr(f\circ g)=0,\;\;
{\rm for\;\;all}\;\;
f\in\c{A}(Y,X)\}.
$$

Morphisms in $\c{N}$ are often called also
{\it numerically trivial}\footnote{
In \cite[pag. 324]{DeSEGL} and \cite[6.1]{DeSt04} they are called ``n\'egligeables''.}.
In the case of motives, indeed, they
correspond to numerically trivial cycles. We refer to \cite[7.1.4, 7.1.5, 7.1.6, 7.4.2.]{AK} 
for other properties of these sets.
Here we just notice that it is not difficult to see the following.

\subsubsection{{\sem Proposition}}\label{NMonidalIdeal}
{\it $\c{N}$ is a monoidal ideal
 of $\c{A}$. If $\c{A}(\I)$ is a field $\c{N}$ is the biggest such ideal, and if moreover
 $\c{A}$ is semisimple then $\c{N}=0$.}

\subsection{Isotypic Schur functors and finiteness conditions}
Assume that $\c{A}(\I)$  contains $\QQ$. 
The partitions $\lambda$ of an integer $|\lambda|=n$ give a complete
set of mutually orthogonal {\it central} idempotents\footnote{Please note that this is not the Young idempotent defining $V_{\lambda}$; it is the {\it central} idempotent defining 
the isotypic component of $V_{\lambda}$ inside the regular representation $\QQ\Sim_{n}$.} 
$
\mathsf{d}_\lambda:=\frac{\dim V_\lambda}{n!}\sum_{\sigma\in \Sigma_n} \chi_\lambda(\sigma)\sigma\in\QQ\Sim_{n}
$
in the group algebra $\Q\Sim_n$ of the symmetric group on $n$ letters with $\QQ$-coefficients (see \cite{FH}), where $\chi_{\lambda}$ is the character of the irreducible representation $V_{\lambda}$ of $\Sigma_{n}$ associated to the partition $\lambda$.
For any $n\in\NN$ and any object $X$ of $\c{A}$, the group $\Sigma_{n}$ acts naturally on $X^{\tens n}$ 
by means of the symmetry of the $\otimes$-category,
we then have also a set of complete mutually orthogonal idempotents indexed by partitions $\lambda$ of $n$
$$
\mathsf{d}_\lambda^{X}:=\frac{\dim V_\lambda}{|\lambda|!}\sum_{\sigma\in \Sigma_n} \chi_\lambda(\sigma)\sigma_{X,\dots,X}\in \c{A}(X^{\tens n})
$$
for each object $X$ of $\c{A}$.
Being $\c{A}$ pseudoabelian, we thus define an endofunctor on $\c{A}$ by setting $S_\lambda(X)=\mathsf{d}^{X}_\lambda(X^{\tens n})$. 
We call it the {\it isotypic Schur functor}; it is a multiple of the classical Schur functor ${\mathbb S}_{\lambda}$ corresponding to $\lambda$ (\cite{FH}, \cite{DeCT02}). 
In particular, we set $\mathrm{Sym}^n(X)=S_{(n)}(X)$ and $\Lambda^n(X)=S_{(1^n)}(X)$.  
Note also that $S_{(0)}=\I$ is the constant functor with value $\I$, and $S_{(1)}=\cid{\c{A}}$.
An object $X$ of $\calA$ is {\it Schur-finite}\footnote{Being $S_{\lambda}(X)={\mathbb S}_{\lambda}(X)^{\oplus \dim V_{\lambda}}$ we have:
$S_\lambda(X)=0\iff {\mathbb S}_{\lambda}(X)=0$.} if there is a partition $\lambda$ such that $S_\lambda(X)=0$.
Schur-finiteness is stable under direct sums, tensor products, duals, and taking direct summands
(see \cite{DeCT02},
\cite{KiFD}, \cite{AK}, \cite{MaSFM}, \cite{DeMaSFN} and \cite{DeMa09} for further reference). 
For further reference we point out the following two propositions.

\subsubsection{{\sem Proposition}}\label{TechnicalProp1}
{\it Let $\lambda$ be a partition of a non negative integer $n$, let $X$ be an object of a $\otimes$-category $\c{A}$,
and let $F\colon \c{A}\To\c{B}$ be a $\otimes$-functor, then:
\begin{enumerate}
\item 
$\chi(S_{\lambda}(X))=
\frac{\dim V_{\lambda}}{|\lambda|!}
\sum_{\sigma\in \Sigma_{|\lambda|}}\chi_{\lambda}(\sigma)\chi(A)^{|\textrm{cycles of }\sigma|}=
\frac{(\dim V_{\lambda})^{2}}{|\lambda|!}{\sf cp}_{\lambda}(\chi(X))$, where $\cp_\lambda(T):=
\prod_{(i,j)\in\lambda}(T+j-i)\in\ZZ[T]$ is the {\rm content polynomial} of the partition $\lambda$.
In particular, if $S_{\lambda}(X)=0$ then $\chi(X)$ is a root of $\cp_\lambda(T)$, hence it is an integer
$\chi(X)\in\{-|\textrm{columns of }\lambda|+1,\dots, |\textrm{rows of }\lambda|-1\}$.
\item $F(S_{\lambda}(X))=S_{\lambda}(F(X))$. 
\end{enumerate}
}
\begin{proof} 
Point (1) follows from part (1) of  Proposition \ref{Abramsky}, which relies on the abstract computation  of \cite[Proposition 3]{Abr}, together with the known properties of the content polynomial of the partition $\lambda$ 
(see \cite[I.1, Example 11, I.3, Example 4 and the proof of I.7(7.6)]{Macd}). 
From the very definitions it's clear that $F(\mathsf{d}^{X}_\lambda)=\mathsf{d}^{F(X)}_\lambda$ hence point (2).
\end{proof}

\subsubsection{{\sem Proposition}}\label{TechnicalProp2}
{\it Let $X=X_{0}|X_{1}$ be a finitely generated supervector space over a characteristic zero field $\FF$, where
$X_{0}$ is the even part and $X_{1}$ is the odd part  of $X$. 
Then $S_{\lambda}(X)=0$ if and only if $\lambda\supseteq ((1+\dim_{\FF} X_{1})^{(1+\dim_{\FF} X_{0})})$, that is
$\lambda$ has at least $\dim_{\FF}(X_{0})+1$ rows and $\dim_{\FF}(X_{1})+1$ columns.
}
\begin{proof} 
It's not difficult to prove it directly, or read \cite[Corollaire 1.9]{DeCT02}.
\end{proof}

\section{Objects (geometrically) of integral type}
In order to keep prerequisites to the minimum I shall not review Kahn's concept of ``multiplicity'', instead I use a definition of being ``(geometrically) of integral type'' which only refers to the Euler characteristics of simple (actually $\varepsilon$-simple) objects.
Moreover, although it is not always necessary, for the pourpous of this note it is enough to work with $\FF$-linear symmetric rigid $\otimes$-categories $\c{A}$ such that $\c{A}(\I)=\FF$ with $\FF$ an {\it algebraically closed} field of characteristic zero (in particular, 
the adverb ``geometrically'' is here pleonastic, see \cite[2.1 d) and 2.2 d)]{KaZFM}).

\subsection{{\bf Definition}} 
\begin{enumerate}
\item
An object $X$ of $\c{A}$ is {\it (geometrically) of integral type} if 
\begin{itemize}
\item[$(i)$] $\c{A}(X)$ is a  finite dimensional semisimple $\FF$-algebra, and
\item[$(ii)$] $\chi(X_{i})=\Tr(\cid{X_{i}})\in\ZZ$ for each 
direct summand $X_{i}$ of $X$ with $\c{A}(X_{i})$ a simple 
$\FF$-algebra.
\end{itemize}
\item
The category $\c{A}$ is said to be  {\it (geometrically) of integral type}
if every object $X$ of $\c{A}$ is such.  The full subcategory
of $\c{A}$ consisting of objects (geometrically) of integral type
is denoted $\c{A}_{\textrm{int}}$.
\end{enumerate}
\subsection{{\bf Remarks}}
\begin{itemize}
\item[(a)] For a review of the notion of a ``semisimple'' category we refer to \cite[2.1.2 and Appendice A]{AK}. 
In \cite[4.2]{Knop}, an object $X$ is called 
{\it $\varepsilon$-semisimple} (resp. {\it $\varepsilon$-simple}) if 
$\c{A}(X)$ is semisimple (resp. simple) ring (without any fineteness assumption). If needed, see also \cite[4.5]{Knop} and 
the comment after \cite[4.7]{Knop} for an explanation of the decomposition of an 
$\varepsilon$-semisimple object into a finite direct sum of 
$\varepsilon$-simple objects.
\item[(b)] 
Note that $\varepsilon$-semisimpleness is quite a mild property, for example if $\c{A}$ is the (not semisimple) rigid $\CC$-linear $\otimes$-category of locally free sheaves (of finite rank) over 
the complex projective line $\PP_{1}$ (or any other geometrically connected scheme with at least two points) then $\I=\c{O}_{\PP_{1}}$ and $\c{A}(\c{O}_{\PP_{1}}(n))=\c{A}(\c{O}_{\PP_{1}})=\c{O}_{\PP_{1}}(\PP_{1})=\CC$ for any $n\in\ZZ$, and the objects 
$\c{O}_{\PP_{1}}(n)$ are of course indecomposable but far from being ``simple'' ({\it i.e.} with no proper subobjects, in the categorical sense): they are not even Artinian!
\item[(c)] 
Being $\c{A}(\I)=\FF$, the (here undefined) ``multiplicity'' $\mu(X)$ of an $\varepsilon$-simple object $X$ 
is exactly the endomorphism $\mu(X)\in\c{A}(X)$
such that $\mu(X)=\chi(X)\cdot \cid{X}$ (see \cite[1.3 b), 2.2 d)]{KaZFM}). 
\end{itemize}

\subsection{{\bf Proposition (Stability properties)}}
{\it  The subcategory $\c{A}_{\textrm{int}}$ is closed under direct sums, direct summands, and duality; moreover it contains the Schur-finite objects. 
}
\begin{proof} The stability properties are proved in \cite[2.2]{KaZFM}.
To see that $\c{A}_{\textrm{int}}$ contains the Schur-finite objects it's enough to note that each direct summand of a Schur-finite
object $X$ is still such and that, using part (1) of Proposition \ref{TechnicalProp1}, the Euler characteristic $\chi(X)$ of a Schur-finite object $X$, say $S_{\lambda}(X)=0$, has to be a root of the content polynomial of the partition $\lambda$, and hence in particular $\chi(X)\in \ZZ$. This is also proved in \cite[Proposition 2.2 e)]{KaZFM} by means of the main, and deep, result of \cite{DeCT02}.
\end{proof}

\subsection{{\bf Two questions}}\label{Qs} The following questions arise naturally:
\begin{itemize}
\item[a)] (Weibel) Does being of integral type imply being Schur-finite?

This question is very interesting in that, although Deligne characterized 
Schur-finite objects in {\it abelian} $\otimes$-categories $\c{A}$, with
$\c{A}(\I)$ an algebraically closed field of characteristic zero, as those whose
lenght of $\otimes$-powers grows at most exponentially (cf. \cite[Proposition 0.5.(i)]{DeCT02}),
there seems to be no analog of \cite[Th\'eor\`eme 7.1]{DeCT90}\footnote{It says the following: {\it In an {\it abelian} $\otimes$-categories $\c{A}$, with
$\c{A}(\I)$ an algebraically closed field of characteristic zero we have: for every $X$ there is an $n\in\NN$ such that $\wedge^{n} X=0$ if and only if $\chi(X)\in\NN$ for every object $X$. This is also equivalent to $\c{A}$ being tannakian.}} for Schur-finiteness.

\item[b)] (Kahn, question to the author) Is $\c{A}_{\textrm{int}}$ closed under $\otimes$ for every $\c{A}$?
\end{itemize}
In this note we answer, in the negative, to these questions as follows.
\begin{itemize}
\item[a)]
In section 3, we show that for each $n\in \ZZ$ there is a {\it not semisimple} $\CC$-linear
rigid, symmetric $\otimes$-categories $\c{T}_n$, 
with hom-sets of finite dimension over $\CC$,
``freely generated'' by an $\varepsilon$-simple 
object $X$ with Euler characteristic 
$\chi(X)=n$. Whence $X$ is of integral type, but it is easy to see that $X$ is not Schur-finite.
This could appear not so definitive in that the $\otimes$-category $\c{T}_n$ in which $X$ sits is not semisimple and not of integral type as a whole\footnote{Note that, by \cite[Proposition 0.5(ii)]{DeCT02} (but it is also easy to see directly), 
any semisimple symmetric rigid $\CC$-linear $\otimes$-category
(without assuming of integral type) having only {\it finitely many simple objects} is necessarily 
Schur-finite, and therefore its pseudoabelian hull (=abelian hull) is of integral type and of {\it homological origin} in the sense
of \cite[Definition 5.1.b)]{KaZFM}
by \cite[Th\'eor\`eme 0.6]{DeCT02}. 
This kind of categories (also not symmetric), called ``fusion categories'', are well known and deeply studied in different contexts (see, e.g., \cite{EtNiOs}), but they are of no use here.

Therefore, if we look for a semisimple $\c{A}$ of integral type, but not Schur-finite,
$\c{A}$ must have infinitely many simple objects,
and, in view of Deligne's Theorem \cite[Proposition 0.5.(i)]{DeCT02}, the length of their tensor powers must have a ``more than exponential'' growth.}. 

In order to give a fully satisfactory example we quickly review Deligne's construction in section 4, and we study, in section 5, a suitable concept of ``\'etale algebras'' in $\otimes$-categories to fully employ its universal property. 
In Proposition \ref{repstnotschur} of section 6 we give two proofs of the fact that the $\otimes$-generator $X=[1]$ of  Deligne's category $\textrm{Rep}(S_{t},\CC)$, defined in \cite{DeSt04}, is not Schur-finite (for all $t\in\CC\setminus \NN$). 
The first proof, in the style of the proof of Proposition \ref{FreeNotSF}, relies on the universal property of $\textrm{Rep}(S_{t},\CC)$:
the key result is Proposition \ref{NoSuper}.
The second proof comes from an explicit lower bound on the rate of growth of $\textrm{length}_{\textrm{Rep}(S_{t},\CC)}([1]^{n})$,
based on two results of Deligne (\cite[Proposition 5.1, Lemme 5.2]{DeSt04}), showing that the length of the tensor powers of $[1]$ grows ``more than exponentially''. The conclusion then follows from \cite[Proposition 0.5.(i)]{DeCT02}.

In Proposition \ref{repstintegraltype} we eventually show that $\textrm{Rep}(S_{t},\CC)$ is (geometrically) of integral type 
if $t\in\ZZ\setminus\NN$\footnote{For $t\in\NN$ the category $\textrm{Rep}(S_{t},\CC)$ is not semisimple.}. 

\item[b)]
In section 7, again working with $\textrm{Rep}(S_{t},\CC)$, 
we show that $\c{A}_{\textrm{int}}$ is not always a $\otimes$-subcategory of $\c{A}$: 
for suitably chosen $t\in\CC\setminus\ZZ$ and simple object $\{\lambda\}\in \textrm{Rep}(S_{t},\CC)$ we show that $\{\lambda\}$ is of integral type but $\{\lambda\}\otimes \{\lambda\}$ is not such. This is achieved by means of Deligne's description of the Grothendieck ring of the semisimple category $\textrm{Rep}(S_{t},\CC)$.
\end{itemize}

\section{Free rigid $\otimes$-categories on one object}

In \cite[Examples (1.26)]{DMOS} there is a construction of the free rigid additive $\otimes$-category on one object $(\c{T}, X_{\c{T}})$.
After recalling what is meant by this we introduce the free rigid $\otimes$-category on one object $\c{T}_{t}$
{\it with prescribed Euler characteristic} $t$ 
and use it to give a first negative answer to Weibel's question \ref{Qs} a).
\medskip

Let $R$ be any commutative ring with identity.

\subsection{{\bf Definition}} A {\it free rigid $R$-linear $\otimes$-category on one object} is a pair $(\c{T}, X_{\c{T}})$ such that
\begin{itemize}
\item[$(i)$]  $\c{T}$ is a rigid $R$-linear $\otimes$-category, and
\item[$(ii)$] for any object $B$ in any rigid $R$-linear $\otimes$-category $\c{B}$
ther is, up to $\otimes$-isomorphism, a unique $\otimes$-functor $F\colon\c{T}\To \c{B}$
such that $F(X_{\c{T}})=B$.
\end{itemize}

\subsection{{\bf Remarks}}
\begin{itemize}

\item[a)]
Such a pair $(\c{T}, X_{\c{T}})$ is clearly unique up to $\otimes$-equivalences.
The existence of such a gadget can be achieved, for example, either by ``abstract non-sense'',
in the vein of \cite{DayNCCC}, or ``constructively'', 
as in \cite[(1.26)]{DMOS}.
\item[b)]
It follows that $X_{\c{T}}$ is a $\otimes$-generator of $\c{T}$ ({\it i.e.}, any object is a direct summand 
of a finite direct sum of tensor powers of $X_{\c{T}}$ and its dual $X_{\c{T}}\dual$) and that $\c{T}(\I)=R[T]$
with $T=\chi(X_{\c{T}})$ algebraically independent over $R$.
\item[c)]
Note that the free $R$-linear $\otimes$-category on one object and the
free {\it rigid} $R$-linear $\otimes$-category $\c{T}$ are related but quite
different: the endomorphisms of $\I$ of the first are reduced to $R$, while
$\c{T}(\I)$ must contain an algebraically independent element over $R$ as remarked above.
\end{itemize}

\subsection{{\bf Definition}} Let $t\in R$. A {\it free rigid $R$-linear $\otimes$-category on one object with Euler characteristic} $t$ is a pair $(\c{T}_{t}, X_{t})$ such that
\begin{itemize}
\item[$(i)$]  $\c{T}_{t}$ is a rigid $R$-linear $\otimes$-category, 
\item[$(ii)$] $\chi(X_{t})=t$, and
\item[$(iii)$] for any object $B$ with $\chi(B)=t$ in any rigid $R$-linear $\otimes$-category $\c{B}$
ther is, up to $\otimes$-isomorphism, a unique $\otimes$-functor $F\colon\c{T}_{t}\To \c{B}$
such that $F(X_{t})=B$.
\end{itemize}
As above, it follows that $X_{t}$ is a $\otimes$-generator of $\c{T}_{t}$ and that $\c{T}_{t}(\I)=R$.

\subsection{{\bf Lemma}}
{\it For any commutative ring with identity and any $t\in R$ there exists 
a free rigid $R$-linear $\otimes$-category on one object with Euler characteristic $t$,
and any two such are $\otimes$-equivalent.}
\begin{proof}
Let $\c{T}$  be the free rigid $R$-linear $\otimes$-category on one object, and let
$\c{I}_t$ be the $\otimes$-ideal $\c{I}_t$ of $\c{T}$ generated by the morphism $T-t\in\c{T}(\I)$, i.e.
$\c{I}_{t}(U,V)=(T-t)\c{T}(U,V)$. 
The pseudo-abelian envelope $\c{T}_{t}$ of
$\c{T}/\c{I}_t$ is a rigid $R$-linear $\otimes$-category 
satifying the required universal property. 
The last part is clear.
\end{proof}

Applying the previous construction to $R=\CC$ we get the first example answering question \ref{Qs} a).

\subsection{{\bf Proposition}}\label{FreeNotSF}
{\it Let $\c{T}_{z}$ to be the free rigid $\CC$-linear $\otimes$-category on one object $X_{z}$ with Euler characteristic $z\in\CC$. Then:
\begin{enumerate}
\item $X_{z}$ is not Schur-finite for any $z\in \CC$.
\item  $X_{z}$ is of integral type if and only if $z\in\ZZ$.
\end{enumerate}
 }
\begin{proof}
(1) If $z\in\CC\setminus\ZZ$ this is obvious, indeed if $S_{\lambda}(X_{z})=0$ for some partition 
$\lambda$ then $z=\chi(X_{z})$ is a root of the content polynomial of $\lambda$ by (1) of Proposition \ref{TechnicalProp1},
and these roots are integers by definition of content polynomial.

In case $z\in\ZZ$, we see that $X_{z}$ cannot be  Schur-finite by universality. 
Assume by contradiction that $S_{\lambda}(X_{z})=0$ for some partition $\lambda$. 
Now, for any $n\in\NN$, by ``definition'' (i.e. the universal property),
there exists a (unique) $\otimes$-functor $F_{n}\colon \c{T}_z\To s\c{V}_{\CC}$ sending the 
$\otimes$-generator $X_{z}$ to the super vector space $F_{n}(X_{z})$ taken to be 
$\CC^{n+z}|\CC^n$ if $z\geq 0$ or $\CC^{n}|\CC^{n+z}$ if $z\leq 0$. 
But by (2) of Proposition \ref{TechnicalProp1} we would also have 
$0=F_{n}(S_{\lambda}(X_{z}))=S_{\lambda}(F_{n}(X_{z}))$.
Taking $n>\max\{|\textrm{rows of } \lambda|, |\textrm{columns of } \lambda|\}$ 
we get in any case a contradiction with Proposition \ref{TechnicalProp2}, 
hence $S_{\lambda}(X_{z})\neq 0$ for any partition $\lambda$.

(2)  $\c{T}_{z}(X_{z})=\CC$ hence $X_{z}$ is $\varepsilon$-simple with $\chi(X_{z})=z$.
\end{proof}

\subsection{{\bf Remark}} Let $z=n\in\ZZ$.
If $n=0$, then $X_{0}$ is a {\it phantom} in $\c{T}_0$, that is $\cid{X_{0}}\in\c{N}$. 
Indeed $\Tr(f)=\Tr(f\cdot\cid{T})=f\cdot\Tr(\cid{T})=0$ for any $f\in\c{T}_{0}(X_{0})=\CC$.
More generally, all the objects of $\c{T}_0$, but those in the subcategory generated by $\I$, are phantoms.
For $n\in\ZZ\setminus\{0\}$, the object $X_{n}$ is {\it not a phantom}.
But all the objects of $\c{T}_n$ represented by a partition (i.e. $S_\lambda(X_{n})$) 
with at least $n+1$ rows if $n>0$, or at least $n+1$
columns if $n<0$, but those in the subcategory generated by $\I$, are phantoms. Indeed these objects are in
the kernel of the $\otimes$-functor sending $X_{n}$ to $\CC^{n}|0$ for $n>0$ and to $0|\CC^{n}$ if $n<0$
(hence the idempotents defining these objects are all universally of trace zero).
Note that, in particular, $X_{n}$ is even or odd in $\c{T}_n/\c{N}_n$ according to $n\geq 0$ or $n\leq 0$.
Hence one cannot deduce a semisimple example from $\c{T}_n$ in this way.

The categories $\c{T}_n$, with $n\in\ZZ$, are clearly {\it neutrally of homological type},
as any $\CC$-super vector space of Euler characteristic $n$ gives a $\otimes$-functor (realization),
but they are not of {\it homological origin} in the sense of \cite[Definition 5.1.]{KaZFM}
(as we see also from $\c{N}_n\neq 0$ for all $n\in\ZZ$).

This ``singular'' behaviour of ``free constructions'' under certain specializations has been already 
explicitely observed and studied by Deligne (cf. \cite{DeSEGL} and \cite{DeSt04}).

\section{Semisimple examples: Deligne's construction $\textrm{Rep}(S_{t},R)$} 
In \cite{DeSt04} (a preprint was available since $2004$) Deligne has constructed several ``interpolating'' families 
of abelian rigid $\otimes$-categories for  categories of representations associated to the  series of 
symmetric, orthogonal, and linear groups. His work has been generilazed by F. Knop (\cite{Knop}).
More specifically, for any commutative ring with identity $R$ and any element $t\in R$, Deligne has constructed a rigid $R$-linear pseudo-abelian $\otimes$-category $\textrm{Rep}(S_{t},R)$
satisfying a suitable universal property. If $R=\FF$ is a field of characteristic zero and $t\in\FF\setminus\NN$ then $\textrm{Rep}(S_{t},R)$ is abelian semisimple, while for $t\in\NN$ then 
$\textrm{Rep}(S_{t},R)/\c{N}$ is $\otimes$-equivalent
to the category of (finite dimensional) $\FF$-linear representations of the symmetric group $\Sigma_{t}$ on $t$ letters.

In this section I will not really enter in either Deligne's nor Knop's construction, 
I shall only  give some hints in the construction
stressing the properties of $\textrm{Rep}(S_{t},R)$ which I will use.

\subsection{Sketch of Deligne's construction $\textrm{Rep}(S_{t},R)$} 

Let $R$ be any commutative ring with identity and let $t\in R$ an element of it.
The rigid $R$-linear pseudo-abelian $\otimes$-cat\-egory $\textrm{Rep}(S_{t},R)$, whose hom-sets are finitely generated projective $R$-module, is obtained  
in three steps:
$$
\textrm{Rep}_{0}(S_{T})\To\textrm{Rep}_{1}(S_{T})\To\textrm{Rep}(S_{t},R),
$$
where $\textrm{Rep}(S_{t},R)$ is defined as the pseudo-abelian envelope of the specialization to $T\mapsto t$ of the
rigid $\ZZ[T]$-linear $\otimes$-category $\textrm{Rep}_{1}(S_{T})$ (\cite[2.16]{DeSt04}),
which is the additive envelope of the $\ZZ[T]$-category $\textrm{Rep}_{0}(S_{T})$ (\cite[2.12]{DeSt04}).

\subsubsection{{\bf The $\ZZ[T]$-category $\textrm{Rep}_{0}(S_{T})$}}\label{Rep0}
The {\it objects} of $\textrm{Rep}_{0}(S_{T})$ are finite sets. 
The symbol $[U]$ denotes the object corresponding to the finite set $U$,
if $U=\{1,\dots, n\}$ one writes $[n]$ for $[U]$.
{\it Morphisms} between  $U$, $V$ (finite sets)
are given by (the $\ZZ[T]$-free module generated by) {\it glueing data} on $U$ and $V$
 (``donn\'ee de recollement sur $U$, $V$''), 
i.e. equivalence relations (equivalently, partitions) $R$ on $U\coprod V$ inducing the discrete equivalence on $U$ and $V$.
In particular the endomorphism ring of the object $[\emptyset]$ is $\textrm{Rep}_{0}(S_{T})([\emptyset])=\ZZ[T]$.
Such morphisms are {\it composed} according to suitable {\it universal polynomial rules} 
(described in \cite[2.10]{DeSt04})  which are products of linear factors of the form 
$T-i$ with $i\in\NN$ (\cite[(2.10.2)]{DeSt04}). 
In this way $\textrm{Rep}_{0}(S_{T})$ is a category enriched over the category of $\ZZ[T]$-modules:
actually its hom-sets $\textrm{Rep}_{0}(S_{T})([U],[V])$ are finitely generated free $\ZZ[T]$-modules 
for any pair of finite sets $U,V$, and composition is $\ZZ[T]$-bilinear. 
Note that $\textrm{Rep}_{0}(S_{T})$ is not an additive category: it lacks products and it is not even pointed
({\it i.e.} there is no zero object).

\subsubsection{{\bf The rigid $\ZZ[T]$-linear $\otimes$-category $\textrm{Rep}_{1}(S_{T})$}}\label{Rep1} 
As remarked above the category $\textrm{Rep}_{0}(S_{T})$ is not even additive, nor monoidal. 
It tourns out that making it additive allows one 
to define a structure of rigid $\otimes$-category  on the resulting 
$\ZZ[T]$-linear category.
The category $\textrm{Rep}_{1}(S_{T})$ is defined as the {\it additive envelope} of $\textrm{Rep}_{0}(S_{T})$, it is therefore
a $\ZZ[T]$-linear category.
Its objects are $n$-tuples of objects of the former category, for $n\in\NN$, 
and morphisms between them are just matrices of morphisms from $\textrm{Rep}_{0}(S_{T})$, 
composed accordingly.
The $\otimes$-structure is defined extending biadditively to  
$\textrm{Rep}_{1}(S_{T})\times \textrm{Rep}_{1}(S_{T})$ the bifunctor
$$
-\otimes-\colon\textrm{Rep}_{0}(S_{T})\times\textrm{Rep}_{0}(S_{T})\To\textrm{Rep}_{1}(S_{T})
$$ 
given by  $[U]\otimes[V]:=\oplus [C]$ where the (formal) direct sum is extended over all
glueing data on $U$ and $V$.
In this way $\I=[\emptyset]$ is the $\otimes$-unit and, for example, $[1]^{\otimes n}$
is the direct sum of $[U/\c{R}]$ over all equivalence relations (equivalently, partitions) 
$\c{R}$ on $U=\{1,\dots, n\}$, in particular $[n]$ is a direct summand of $[1]^{\otimes n}$.

A closer inspection to the definitions (cf. \cite[2.16]{DeSt04}, or \cite[(3.15)]{Knop}) shows that each object $[U]$ of $\textrm{Rep}_{0}(S_{T})$ (and hence of $\textrm{Rep}_{1}(S_{T})$)
is actually {\it selfdual}: $[U]\dual=[U]$. 
Hence the object $[1]$ is a $\otimes$-generator of $\textrm{Rep}_{1}(S_{T})$.

\subsubsection{{\bf The rigid $R$-linear pseudo-abelian $\otimes$-category $\textrm{Rep}(S_{t},R)$}} 
The ca\-tegory $\textrm{Rep}(S_{t},R)$ is defined as the pseudo-abelian envelope of the {\it specialization}
$\textrm{Rep}_{1}(S_T)\otimes_{T\mapsto t}R$. If $U$ is a finite set its corresponding object in $\textrm{Rep}(S_{t},R)$ is $[U]_{t}$, and it will be denoted simply $[U]$ if $t$ is clear from the context.
By construction it is a rigid $R$-linear pseudoabelian $\otimes$-category, 
and it is not difficult to check that its hom-sets are finitely generated projective $R$-modules. 
As for $\textrm{Rep}_{1}(S_T)$
the unit object of $\textrm{Rep}(S_{t},R)$ is $\I=[\emptyset]_{t}$, with $\textrm{Rep}(S_{t},R)(\I)=R$,
and every object $[U]$ ($U$ a finite set)
is self dual and endowed with a natural structure of ``ACU'' algebra\footnote{ACU stands for ``associatif, commutatif, \`a unit\'e''
as in \cite{SaaCT} and \cite[1.9]{DeSt04}} ({\it i.e.} commutative monoid) in
$\textrm{Rep}(S_{t},R)$ (see \cite[1.2, 2.5, 2.16]{DeSt04}).
The object $[1]$ is a self dual algebra  with $\chi([1])=t$ such that $\textrm{Rep}(S_{t},R)$ is the 
pseudoabelian hull of the full subcategory on the objects $[1]^{\otimes n}$ (with $n\in\NN$), hence $[1]$ is still a $\otimes$-generator of
$\textrm{Rep}(S_{t},R)$.

\subsubsection{{\bf  ``Partition objects'' and their Euler characteristic}}\label{SimpleObjects} 
Let $\lambda=(\lambda_{1}\geq \lambda_{2}\geq\dots\geq 0)$
be a partition (of some integer $|\lambda|=\sum_{i}\lambda_{i}$). 
Assume that $R$ is a $\QQ$-algebra and that $t\in R$ is such that $t-n\in R$ is a unit
for any $n\in\NN\cap [0,2|\lambda|-2]$, then by \cite[Remarque 5.6]{DeSt04} there is, functorially w.r.t. $R$, 
an object $[\lambda]_{t}$ of $\c{A}=\textrm{Rep}(S_{t},R)$ such that: 
\begin{itemize}
\item[a)]
$\c{A}([\lambda]_{t})=R$, 
\item[b)]
if $\lambda,\mu$ are two distinct partitions, and  $t-n\in R$ is a unit
for any $n\in\NN\cap [0,2(|\lambda|\vee|\mu|)-2]$, then
$\c{A}([\lambda]_{t}, [\mu]_{t})=0$,
\item[c)] if $t-n\in R$ is a unit
for any $n\in\NN$, then
for every object $Y$ of $\c{A}$ we have:
$$
\oplus_{\lambda} \c{A}([\lambda]_{t}, Y)\otimes_{R}
[\lambda]_{t}\cong Y,
$$
where $\c{A}([\lambda]_{t}, Y)\otimes_{R}
[\lambda]_{t}$ is the object of $\c{A}$
representing the functor (cf. \cite[3.7]{DeSt04})
$$
Z\mapsto \Hom_{R}(\c{A}([\lambda]_{t}, Y),
\c{A}([\lambda]_{t},Z)).
$$
If $R$ is a PID then $\c{A}([\lambda]_{t}, Y)\otimes_{R}
[\lambda]_{t}$ is just a sum of copies of $[\lambda]_{t}$.
\end{itemize}
For example, taking $R=\CC$ and $t\in\CC\setminus\NN$ one has: $\I=[\emptyset]=[(0)]_{t}$, while $[1]=\I\oplus
[(1)]_{t}$ (see \cite[5.1, 5.5]{DeSt04}).

Moreover, by \cite[Lemme 7.3, 7.4]{DeSt04}, for any partition $\lambda$ there is a universal polynomial 
$$
Q_\lambda(T):=\frac{\dim V_\lambda}{|\lambda|!}\prod_{a=1}^{|\lambda|}(T-(|\lambda|+\lambda_a-a))\in\QQ[T]
$$
such that the Euler characteristic of the object $[\lambda]_{t}$ of $\textrm{Rep}(S_{t},R)$ corresponding to $\lambda$ is given by
$
\chi([\lambda]_{t})=Q_\lambda(t)\;\textrm{ for any }t\in R\setminus\NN\cap[0,2|\lambda|-2].
$

\subsection{Theorem (Deligne)}\label{DeligneSemisimple}
{\it
If $R=\FF$ is a field of characteristic zero then
$\textrm{Rep}(S_{t},\FF)$ is a semisimple category (hence also abelian)
if and only if $t\in\FF\setminus\NN$; the simple objects of $\textrm{Rep}(S_{t},\FF)$
are given by partitions and they are {\it absolutely simple}.}
\begin{proof} This is \cite[Th\'eor\`eme 2.18, Proposition 5.1, Th\'eor\`eme 6.2]{DeSt04}. The simple objects are those described in \ref{SimpleObjects}, which are absolutely simple as
$\textrm{Rep}(S_{t},\FF)([\lambda]_{t})=\FF$ 
for any characteristic zero field $\FF$, any partition $\lambda$ and any $t\in\FF\setminus\NN$.
\end{proof}

\section{{\bf Universal property of $\textrm{Rep}(S_{t},R)$: \'etale algebras in rigid $\otimes$-categories}}
We would like to show that $\textrm{Rep}(S_{t},\CC)$ is not Schur-finite for $t\in\CC\setminus\NN$. This is clear if $t\in\CC\setminus \ZZ$, for $[1]$ Schur-finite would implies $t=\chi([1])\in\ZZ$, as already remarked in the proof of Proposition 2.3. 
It remains the case $t\in\ZZ\setminus\NN$.  A way to show that $\textrm{Rep}(S_{t},\CC)$ is not Schur-finite in this case, in analogy with the strategy of section 3, is by means of its universal property which says what are the possibile $\otimes$-functors from $\textrm{Rep}(S_{t},R)$. It is also possible to avoid the use of the universal property, as we will see in the two proofs of Proposition \ref{repstnotschur}, nonetheless we found of some interest the following results.

\subsection{The universal property of $\textrm{Rep}(S_{t},R)$}
Following \cite[8, 8.1, 8.2, 8.3]{DeSt04}
the category $\textrm{Rep}(S_{t},R)$ is characterized by the following universal property:
$R$-linear $\otimes$-functors $F\colon \textrm{Rep}(S_{t},R)\To \c{A}$
from it to any other $R$-linear pseudo-abelian $\otimes$-category $\c{A}$ correspond bijectively, via $A=F([1])$,
with the ``ACU'' algebras ({\it i.e.}, commutative monoids) $(A,m,u)$ of $\c{A}$ such that $\chi(A)=t$ and
$A\otimes A\xrightarrow{m} A\xrightarrow{T_m} \I$ is the self duality of $A$ where $T_m\colon A\To \I$ is the composite\footnote{
More generally, any morphism $m\colon A\otimes M\To M$ 
(think of $M$ as an ``$A$-module'' as in \cite[18]{MacLCA}) in a rigid $\c{A}$ induces
a ``trace'' 
$$
T_{m}:=(A=A\otimes \I\xrightarrow{A\otimes \eta_{M}} A\otimes M\otimes M\dual
\xrightarrow{m\otimes M\dual} M\otimes M\dual
\xrightarrow{\tau_{M,M\dual}}M\dual\otimes M
\xrightarrow{\varepsilon_{M}}\I).
$$
}
$$
A=A\otimes \I\xrightarrow{A\otimes \eta_{A}} A\otimes A\otimes A\dual
\xrightarrow{m\otimes A\dual} A\otimes A\dual
\xrightarrow{\tau_{A,A\dual}}A\dual\otimes A
\xrightarrow{\varepsilon_{A}}\I.
$$
\subsubsection{{\bf Remark}}\label{etalealgebras}
From a categorical perspective, the algebras as above are all obviously (commutative) {\it Frobenius algebras} (\cite{StFM},\cite{FuStin}, \cite{LamLMR}), but the kind of algebras one can get as $F([1])$ is even more restricted.

For example, in the rigid $\FF$-linear $\otimes$-category $\c{A}=\c{V}_{\FF}$  of (finitely generated) vector spaces over the field $\FF$ it is easy to see that the monoid objects $(A, m, u)$ with the properties stated above are exactly the {\it \'etale} algebras over $\FF$ (see \cite[V, p.48, Prop. 1]{BouAII}), equivalently $A\cong \KK_{1}\times\cdots\times\KK_{n}$
where $\KK_{i}$ is any separable finite extension of the field $\FF$.
Note in particular that such an $A$ has to be a reduced ring ({\it i.e.}, it has no non zero nilpotent elements).
\medskip

\subsection{ \'Etale algebras in rigid $\otimes$-categories}
In order to further investigate the kind of algebras on can get as $F([1])$
in a general rigid $R$-linear $\otimes$-category $\c{A}$, with
$F\colon \textrm{Rep}(S_{t},R)\To\c{A}$ a $\otimes$-functor, 
let us fix some terminology allowing us to do some ``multilinear algebra'' in a $\otimes$-category.

\subsubsection{{\bf Definition}}\label{definitionetalealgebra} 
Let $\c{A}$ be a $\otimes$-category.
A {\it bilinear pairing} in a $\otimes$-category $\c{A}$ is any morphism of the form $b\colon X\otimes Y\To \I$. We say that $b$ is a {\it bilinear form} if $Y=X$. If $\c{A}$ is symmetric, we say that a bilinear form $b\colon X\otimes X\To \I$ is a {\it symmetric} [{\it antisymmetric}] if $b\circ \tau_{X,X}=b$ [$b\circ \tau_{X,X}=-b$].
If $\c{A}$ is rigid, we say that $b\colon X\otimes Y\To \I$ is {\it non degenerate}, or a {\it perfect pairing}, if the induced morphism
$$
X=X\otimes\I\xrightarrow{X\otimes \eta_{Y}} X\otimes Y\otimes Y\dual \xrightarrow{b\otimes Y\dual} \I\otimes Y\dual=Y\dual
$$
is an isomorphism.

\subsubsection{{\bf Remark}} The reader will have no difficulties in cheking that in the rigid abelian
$\otimes$-category $\c{A}=\c{V}_{\FF}$  of (finitely generated) vector spaces over the field $\FF$ all the previous notions coincide with the usual ones. 

\subsubsection{{\bf Definition}} An {\it \'etale} algebra in a rigid $\otimes$-category $\c{A}$ is a monoid object $(A,m,u)$ in $\c{A}$
such that $m$ is commutative and the symmetric bilinear form 
$$
A\otimes A\xrightarrow{m} A\xrightarrow{T_m} \I,
$$ 
where $T_{m}\colon A\To \I$ is the composite
$$
A=A\otimes \I\xrightarrow{A\otimes \eta_{A}} A\otimes A\otimes A\dual
\xrightarrow{m\otimes A\dual} A\otimes A\dual
\xrightarrow{\tau_{A,A\dual}}A\dual\otimes A
\xrightarrow{\varepsilon_{A}}\I,
$$
is non degenerate, {\it i.e.} 
$$
((T_{m}\circ m)\otimes\cid{A\dual})\circ (\cid{A}\otimes \eta_{A})\colon 
A=A\otimes \I\xrightarrow{A\otimes \eta_{A}}A\otimes A\otimes A\dual
\xrightarrow{(T_{m}\circ m)\otimes A\dual} A\dual
$$
is an isomorphism $A\cong A\dual$.

\subsubsection{{\bf Remark}} 
We want to stress, following \cite[2]{JSV} and \cite[4.3]{Abr} that given any $m\colon A\otimes M\To M$ its ``trace'',
as defined above, is nothing but
$$
T_{m}=\Tr^{M}_{A, \I}(m\colon A\otimes M\To M\cong\I\otimes M).
$$
For example, if $A$ is an algebra of finite dimension over a field and $M$ is a finitely generated $A$-module 
then $T_{m}=\chi_{M}$ is the usual character of the module $M$.
\medskip

We already mentioned that in the rigid $\otimes$-category $\c{V}_{\FF}$ of finitely generated vector spaces
over $\FF$ an \'etale algebra in this categorical sense is nothing but an \'etale algebra in the classical sense.
Let us see what are the \'etale algebras in the rigid $\otimes$-category $\sf{s}\!\c{V}_{\FF}$ of super vector spaces.

\subsection{\'Etale algebras in supercategories}
The {\it super category} of a $\otimes$-category $\c{A}$ is the $\otimes$-category $\sf{s}\!\c{A}$ whose
objects and morphisms are as in the product category $\c{A}\times\c{A}$, but the {\it tensor product}
is defined on objects as
$$
X\otimes_{\sf{s}} Y:=
X_{0}\otimes Y_{0}\,\oplus\, X_{1}\otimes Y_{1} 
\mid
X_{0}\otimes Y_{1}\,\oplus\, X_{1}\otimes Y_{0}
$$
if $X=X_{0}\mid X_{1}$ and $Y=Y_{0}\mid Y_{1}$; 
on morphisms $f=f_{0}\mid f_{1}\colon X\To Y$ and  $g=g_{0}\mid g_{1}\colon A\To B$ we have
$$
f\otimes_{\sf{s}}g:=
\begin{pmatrix}
f_{0}\otimes g_{0}&0\\
0& f_{1}\otimes g_{1}
\end{pmatrix} 
\Big\arrowvert
\begin{pmatrix}
f_{0}\otimes g_{1}& 0\\
0&f_{1}\otimes g_{0}
\end{pmatrix}
\colon X\otimes_{\sf{s}} A\To Y\otimes_{\sf{s}} B.
$$

In this way $\I_{\sf{s}}=\I \mid 0$ and its endomorphism ring is $\sf{s}\!\c{A}(\I_{\sf{s}})=\c{A}(\I)$.

\subsubsection{{\bf Writing convention}}
In what follows it will be helpful to adopt the usual convention about
morphisms between finite biproducts in an additive category:
$$
( X_{1}\oplus\cdots\oplus X_{n}\stackrel{F}{\To} Y_{1}\oplus\cdots\oplus Y_{m})\mapsto
\begin{pmatrix}
F^{1}_{1}&\cdots& F^{n}_{m}\\
\vdots&\ddots&\vdots\\
F^{1}_{m}&\cdots&F^{n}_{m}
\end{pmatrix},
\quad\textrm{with}\quad X_{i}\stackrel{F^{i}_{j}}{\To}  Y_{j}.
$$
Note the following special cases in $\sf{s}\!\c{A}$:
$$
(A\otimes_{\sf{s}} X\stackrel{f}{\To} B\otimes_{\sf{s}} Y)\mapsto 
\begin{pmatrix}
f^{00}_{00}&f^{11}_{00}\\
f^{00}_{11}&f^{11}_{11}
\end{pmatrix}
\Big\arrowvert
\begin{pmatrix}
f^{01}_{01}&f^{10}_{01}\\
f^{01}_{10}&f^{10}_{10}
\end{pmatrix}
\quad\textrm{with}\quad A_{i}\otimes X_{j}\xrightarrow{f^{ij}_{hk}}B_{h}\otimes Y_{k},
$$
$$
(A\otimes_{\sf{s}} B\stackrel{m}{\To} C)\mapsto 
\begin{pmatrix}
m^{00}_{0}&m^{11}_{0}
\end{pmatrix}
\Big\arrowvert
\begin{pmatrix}
m^{01}_{1}&m^{10}_{1}
\end{pmatrix}
\quad\textrm{with}\quad A_{i}\otimes B_{j}\xrightarrow{m^{ij}_{k}} C_{k},
$$
and
$$
(A\stackrel{n}{\To} B\otimes_{\sf{s}} C)\mapsto 
\begin{pmatrix}
n^{0}_{00}\\
n^{0}_{11}
\end{pmatrix}
\Big\arrowvert
\begin{pmatrix}
n^{1}_{01}\\
n^{1}_{10}
\end{pmatrix}
\quad\textrm{with}\quad A_{i}\xrightarrow{n^{i}_{jk}} B_{j}\otimes  C_{k},
$$

\subsubsection{{\bf Associativity}}
It's worth to point out explicitly the structure of the ternary associativity isomorphisms 
$$
\alpha^{\sf{s}}\colon (X\otimes_{\sf{s}} Y)\otimes_{\sf{s}} Z\To
X\otimes_{\sf{s}} (Y\otimes_{\sf{s}} Z).
$$
Writing only the parity symbols,  the structure of $\alpha^{\sf{s}}$ is as follows
$$
\frac{(\alpha^{\sf{s}})_{0}}{(\alpha^{\sf{s}})_{1}}\colon 
\frac{(00)0+(11)0+(01)1+(10)1}{(00)1+(11)1+(01)0+(10)0}
\To
\frac{0(00)+0(11)+1(01)+1(10)}{0(01)+0(10)+1(00)+1(11)} 
$$
Hence $(\alpha^{\sf{s}})_{i}\colon [(X\otimes_{\sf{s}} Y)\otimes_{\sf{s}} Z]_{i}\To
[X\otimes_{\sf{s}} (Y\otimes_{\sf{s}} Z)]_{i}
$ is obtained by ``composing'' the associativity isomorphisms of $\c{A}$
with the cyclic permutation $(243)$ whose matrix representation is
$$
P=\begin{pmatrix}
1&0&0&0\\
0&0&1&0\\
0&0&0&1\\
0&1&0&0
\end{pmatrix}
\quad\textrm{with inverse}\quad
P^{-1}={}^{t}\!P=\begin{pmatrix}
1&0&0&0\\
0&0&0&1\\
0&1&0&0\\
0&0&1&0
\end{pmatrix}.
$$
Note that writing down the matrix representation of a morphisms of the form
$$f\colon A\otimes_{\sf{s}} B\otimes_{\sf{s}} C\To X\otimes_{\sf{s}} Y\otimes_{\sf{s}} Z$$
one should take into account the choice of the brackets.

\subsubsection{{\bf Simmetry}}
It is crucial to recall that the {\it simmetry} of $\sf{s}\!\c{A}$ is defined by the Koszul rule
$$
\tau^{\sf{s}}_{X,Y}=
\begin{pmatrix}
\tau_{X_{0},Y_{0}}&0\\
0& -\tau_{X_{1},Y_{1}}
\end{pmatrix} 
\Big\arrowvert
\begin{pmatrix}
0& \tau_{X_{1},Y_{0}}\\
\tau_{X_{0},Y_{1}}&0
\end{pmatrix}
\colon 
X\otimes Y\cong Y\otimes X.
$$
\subsubsection{{\bf Rigidity}} $\c{A}$ is rigid if and only if $\sf{s}\!\c{A}$ is such and
$(X_{0}\mid X_{1})\dual=X_{0}\dual\mid X_{1}\dual$ 
with evaluation and coevaluation given by
$$
\varepsilon^{\sf{s}}_{X}=
\begin{pmatrix}
\varepsilon_{X_{0}} & \varepsilon_{X_{1}}
\end{pmatrix}
\Big\arrowvert 
0
\colon X\dual\otimes X \To\I,
\quad
\eta^{\sf{s}}_{X}=
\begin{pmatrix}
\eta_{X_{0}} \\ 
\eta_{X_{1}}
\end{pmatrix}
\Big\arrowvert 
0
\colon \I\To X\otimes X\dual.
$$

\subsubsection{{\bf Remark}}\label{tracecommutativemonoid}
Let $m\colon A\otimes_{\sf{s}} M\To M$ be a morphism in $\sf{s}\!\c{A}$ described as
$$
m=
\begin{pmatrix}
m^{00}_{0} & m^{11}_{0}
\end{pmatrix}
\Big\arrowvert
\begin{pmatrix}
m^{01}_{1} & m^{10}_{1}
\end{pmatrix}
\quad\textrm{with}\quad
m^{ij}_{k}\colon A_{i}\otimes M_{j}\To M_{k}.
$$
The trace induced by $m\colon A\otimes_{\sf{s}} M\To M\cong\I_{\sf{s}}\otimes_{\sf{s}} M$ is 
$$
{}^{\sf{s}}T_{m}={}^{\sf{s}}\Tr^{M}_{A,\I}(m)=
\Tr^{M_{0}}_{A_{0},\I}(m^{00}_{0})-\Tr^{M_{1}}_{A_{0},\I}(m^{01}_{1})\mid 0=
T_{m^{00}_{0}}-T_{m^{01}_{1}}\mid 0\colon A_{0}\mid A_{1}\To \I\mid 0.
$$

In case $M=A$ then $m\colon A\otimes_{\sf{s}}A\To A$ and requiring such an $m$ to be 
{\it commutative in $\sf{s}\!\c{A}$} means $m\circ \tau^{\sf{s}\!\c{A}}_{A,A}=m$.
Since 
$$
m\circ \tau^{\sf{s}\!\c{A}}_{A,A}=
\begin{pmatrix}
m^{00}_{0}\circ \tau_{A_{0},A_{0}} & 
-m^{11}_{0}\circ \tau_{A_{1},A_{1}}
\end{pmatrix}
\Big\arrowvert
\begin{pmatrix} 
m^{10}_{1}\circ \tau_{A_{0},A_{1}}&
m^{01}_{1}\circ \tau_{A_{1},A_{0}} 
\end{pmatrix}
$$
this last condition in $\sf{s}\!\c{A}$  is equivalent to the following conditions in $\c{A}$:
$m^{00}_{0}$ is commutative, $m^{11}_{0}$ is anticommutative and 
$m^{01}_{1}\circ \tau^{\c{A}}_{A_{1},A_{0}}=m^{10}_{1}$. 

\subsubsection{{\bf Proposition}}\label{NoSuper}
{\it The only \'etale algebras $A$ in $\sf{s}\c{V}_{\FF}$ are the \'etale algebras of $\c{V}_{\FF}$ thought of as (purely) even superobjects, that is in particular $A_{1}=0$, and hence $\chi(A)>0$.
In particular, if $t\in\FF\setminus\NN$ there are no $\otimes$-functor $F\colon \textrm{Rep}(S_{t},\FF)\To \sf{s}\c{V}_{\FF}$.}
\begin{proof} 
Let $\c{A}=\c{V}_{\FF}$ be the rigid  $\otimes$-category of (finitely generated) vector spaces over the field $\FF$, 
and let $A$ be an  \'etale algebra in $\sf{s}\!\c{A}$.
then  $\Tr^{\sf{s}}_m\circ m\colon A\otimes_{\sf{s}} A\To \I_{\sf{s}}=\FF\mid 0$ in $\sf{s}\!\c{A}$
induces not degenerate symmetric bilinear forms
$(\Tr_{m_{0}^{00}}-\Tr_{m_{1}^{01}})\circ m_{0}^{ii}\colon A_{i}\otimes A_{i}\To \FF$ in $\c{A}$, 
where $\Tr_{m^{0i}_{i}}\colon A_{0}\To \FF$ is the trace map induced by $m^{0i}_{i}\colon A_{0}\otimes A_{i}\To A_{i}$ 
But, by supercommutativity, the image of $m^{11}_{0}\colon A_{1}\otimes A_{1}\To A_{0}$
would be made of elements with square zero, having therefore zero traces. Hence $A_{1}=0$.

Let us assume now by contradiction that there is a $\otimes$-functor 
$F\colon \textrm{Rep}(S_{t},\FF)\To \sf{s}\c{V}_{\FF}$.
Then $F([1])$ would be an \'etale superalgebra of superdimension $\chi(F([1]))=\chi([1])=t<0$,  
{\it i.e.} $F([1])_{1}\neq 0$ which is impossible.
\end{proof}

\section{$\textrm{Rep}(S_{t},\CC)$ is an abelian semisimple not Schur-finite category of integral type for any $t\in\ZZ\setminus\NN$}\label{NotSF}

\subsection{}
We show, with two different proofs, that the $\otimes$-generator $[1]$ of $\textrm{Rep}(S_{t},\FF)$ is not Schur-finite
for any $t\in\FF\setminus\NN$. Note that
in case $t\not\in\ZZ$ this follows easily from part (1) of Proposition \ref{TechnicalProp1} for $t=\chi([1])$.

\subsubsection{{\bf Proposition}}\label{repstnotschur}
{\it Let $\FF$ be an algebraically closed field of characteristic zero, and let $t\in\FF$.
If $t\not\in\NN$ then the object $[1]$ is not Schur-finite in $\textrm{Rep}(S_{t},\FF)$.
}

\subsubsection*{ Proof 1}
Let $t\in\FF\setminus\NN$, if $[1]$ be Schur-finite then by Theorem \cite[0.6]{DeCT02} there is a $\otimes$-functor 
$
F\colon \textrm{Rep}(S_{t},\FF)\To \sf{s}\c{V}_{\FF}
$, 
but this is impossible by Corollary \ref{NoSuper}.
\qed

\subsubsection*{ Proof 2}
As promised, we give also an alternate proof of the fact that $[1]$  
is not Schur-finite in $\textrm{Rep}(S_{t},\FF)$ with $t\in\FF\setminus \NN$ without any reference to
the universal property. 
By \cite[Th\'eor\`eme 2.18]{DeSt04} the category $\c{A}:=\textrm{Rep}(S_{t},\FF)$
is abelian semisimple with simple objects which are absolutely simple, hence for any object $A$ of $\c{A}$ 
we have the lower bound:
$
\textrm{length}_{\c{A}}(A)\geq \sqrt{\dim_{\FF}\c{A}(A)}.
$
Let now $A=[1]^{n}$ with $n\in \NN$. By the very definition of the $\otimes$-product of $\textrm{Rep}_{1}(S_{T})$,
the object $[1]^{\otimes n}$ has $[n]$ as a direct summand.
Moreover, by Deligne's \cite[Proposition 5.1]{DeSt04}, inside the object $[n]$
there is another direct summand: the object $[n]^{*}$, which has the property that 
$\c{A}([n]^{*})\cong \FF\Sigma_{n}$ as algebras, as proved in \cite[Lemme 5.2]{DeSt04}. 
Therefore we have:
$$
\textrm{length}_{\c{A}}([1]^{\otimes n})\geq 
\textrm{length}_{\c{A}}([n])\geq 
\textrm{length}_{\c{A}}([n]^{*})\geq \sqrt{\dim_{\FF}\c{A}([n]^{*})}=
\sqrt{n!}
$$
hence $[1]$ can't be Schur-finite because $\textrm{length}_{\c{A}}([1]^{\otimes n})$ grows faster than exponentially 
(cf. \cite[Proposition 0.5. (i)]{DeCT02}) in the category $\c{A}$ which is a ``cat\'egorie $\FF$-tensorielle'' in the sense of \cite{DeCT02}.
\qed

\subsubsection{{\bf Remark}} We note that
this disproves the claims \cite[4.4, 4.5]{HoHai}. 

\subsection{}
It remains to show that 
$\textrm{Rep}(S_{t},\FF)$ is (an abelian semisimple category) of {\it integral type} for any field of characteristic zero $\FF$ and any $t\in\FF$ such that $t\in\ZZ\setminus\NN$.
As recalled above (see \ref{SimpleObjects} and Theorem \ref{DeligneSemisimple}) the simple objects of $\textrm{Rep}(S_{t},\FF)$ are   given by partitions and they are absolutely simple, hence to show they are geometrically of integral type  it's enough to show they have integer Euler characteristic.

\subsubsection{{\bf Proposition}}\label{repstintegraltype}
{\it Let $\FF$ be an algebraically closed field of characteristic zero, and let $t\in\FF\setminus\NN$. 
Then $\textrm{Rep}(S_{t},\FF)$ is of integral type if and only if $t\in \ZZ\setminus\NN$.}
\begin{proof}
As recalled in \ref{SimpleObjects},
for the simple object $[\lambda]_{t}$ associated to the partition $\lambda$,
Deligne has given the explicit expression 
$$
\chi([\lambda]_{t})=
Q_\lambda(t)\;\textrm{ where }\;
Q_\lambda(T)=
\frac{\dim V_\lambda}{|\lambda|!}\prod_{a=1}^{|\lambda|}(T-(|\lambda|+\lambda_a-a))\in\QQ[T].
$$
We claim that $Q_{\lambda}(\ZZ)\subset\ZZ$.
By \cite[6.4, (7.4.1)]{DeSt04}, $Q_{\lambda}(n)
=\dim V_{\{\lambda\}_n}\in\NN$ for all $n\geq 2|\lambda|+1$,
where $\{\lambda\}_{m}$ is a partition attached to $\lambda$ for each 
$m\geq |\lambda|+\lambda_{1}$ 
(\cite[(6.3.1)]{DeSt04}).
Then the claim follows by the elementary fact recalled in Lemma \ref{integerpolynomial} (I am certain it is well known,
but I do not know any reference).
Therefore all simple objects of $\textrm{Rep}(S_{t},\FF)$
have integer Euler characteristic and $\textrm{Rep}(S_{t},\FF)$ is geometrically of integral type.
\end{proof}

\subsubsection{{\bf Lemma}}\label{integerpolynomial}
{\it
Let $p(T)\in\QQ[T]$ be a polynomial of degree $d$. 
Then $p(\ZZ)\subseteq \ZZ$ if and only if
$p$ takes integer values on $d+1$ consecutive integer points.}
\begin{proof}
By induction on the degree $d$ of $p$. The assertion is true if $d=0$.
Let $d>0$ and assume the result on all polynomials $q(T)\in\QQ[T]$ of degree $d-1$. 
By hypothesis there is an $n\in\ZZ$ such that $p(n+i)\in\ZZ$ for $i\in\{0,\dots, d\}$.
Then, to low the degree, note that $\Delta p(T):=p(T)-p(T-1)\in\QQ[T]$
is a polynomial of degree $d-1$ such that: 
$p(\ZZ)\subseteq \ZZ$ if and only if
$p$ takes an integer value in at least one integer and $\Delta p(\ZZ)\subseteq \ZZ$.
Indeed clearly $p(\ZZ)\subseteq \ZZ$ implies $\Delta p(\ZZ)\subseteq \ZZ$; conversely,
assume $\Delta p(\ZZ)\subseteq \ZZ$ and that there is a $n\in\ZZ$
such that $p(n)\in \ZZ$. Then $p(n+i)=p(n)+\sum_{j=1}^i \Delta p(n+j)\in\ZZ$
for each $i\in \NN_{+}$, and 
$p(n-i)=p(n)-\sum_{j=0}^i\Delta p(n-j)\in\ZZ$ for each $i\in \NN$.
Whence $p(\ZZ)\subseteq \ZZ$.
We conclude the induction as  follows: 
by hypothesis on $p$ we have $p(n)\in\ZZ$ and $\Delta p(T)$ assumes integer values on $[n+1,\dots, n+d]\cap \ZZ$,
hence $\Delta p(\ZZ)\subseteq \ZZ$ by the inductive hyphothesis.
\end{proof}

\section{Tensor powers of objects of integral type need not be such} 

The tensor powers of a geometrically integral (even absolutely simple) object
are not, in general, still such. We can look for an example as follows.

We know that $\textrm{Rep}(S_t, \CC)$, with $t\in\CC\setminus\NN$, is an abelian semisimple rigid $\CC$-linear category, 
and we also have a complete description of its (absolutely) simple objects and their Euler characteristic. 
By \cite[5.10]{DeSt04}, simple objects $[\mu]_{t}$ are given by partitions $\mu$, 
and  $\chi([\mu]_{t})=\frac{\dim V_\mu}{|\mu|}\prod_{a=1}^{|\mu|}(t-|\mu|-\mu_a+a)$. 
Moreover, the Grothendieck ring of such a category is also described in \cite{DeSt04} 
by means of classical representation theory of symmetric groups, 
in particular we can effectively check if a simple summand $[\mu]_{t}$ appears in objects like 
$[\lambda]_{t}\otimes[\nu]_{t}$.

Hence, the idea is to take an absolutely simple object $[\lambda]_{t}$ in a suitable specialization
of $\textrm{Rep}_{1}(S_T)$, with $T\mapsto t\in\CC\setminus\ZZ$,
 in such a way that:
 \begin{itemize}
 \item[(a)]
 $\chi([\lambda]_{t})\in \ZZ$, so that our 
simple object $[\lambda]_{t}$ is geometrically of integral type, but

 \item[(b)] $[\lambda]_{t}\otimes[\lambda]_{t}$
is not geometrically of integral type, i.e. some simple object $[\mu]_{t}$ 
with non integer Euler characteristic appears in its decomposition.
 \end{itemize}
 
Let us work out the details.
 \begin{itemize}
 \item[(a)] Take $\lambda:=(2,1)$. 
Then $\chi([\lambda]_{t})=\frac{1}{3}(t-4)(t-2)t$ for any $ t\in\CC\setminus\ZZ$. 
As $\chi([\lambda]_{-1})=-5$, 
then $\chi([\lambda]_{\tau})=-5$ also for $\tau\in\{\frac{7+i\sqrt{11}}{2},\frac{7-i\sqrt{11}}{2}\}$, i.e.
the other roots of the polynomial $\frac{1}{3}(t-4)(t-2)t+5=\frac{1}{3}(t+1)(t^2-7t+15)$. 
Let me fix $\tau:=\frac{7+i\sqrt{11}}{2}$, from now on I work in the
abelian semisimple rigid $\CC$-linear category $\c{A}:=\textrm{Rep}(S_{\frac{7+i\sqrt{11}}{2}},\CC)$.
Hence $[\lambda]_{\tau}$ is a geometrically integral
object of $\c{A}$.
 \item[(b)] Consider the second tensor power $[\lambda]_{\tau}\otimes[\lambda]_{\tau}$ of $[\lambda]_{\tau}$.
By \cite[5.11]{DeSt04} and the known properties of Littelwood-Richardson coefficients, 
we know that among the simple objects $[\mu]_{\tau}$ appearing in the decomposition of 
$[\lambda]_{\tau}\otimes[\lambda]_{\tau}$  there are some $[\mu]_{\tau}$ with $|\mu|= 6$. 
Hence it is enough to check that 
$\chi([\mu]_{\tau})\not\in\ZZ$ for at least one such $[\mu]_{\tau}$, 
which is indeed the case: take $\mu:=(3,2,1)$ then 
$\left[\mathrm{Ind}^{\Sigma_{6}}_{\Sigma_{3}\times \Sigma_{3}}(V_\lambda\otimes V_\lambda)\colon V_\mu\right]=2$, so the simple object $[\mu]_{\tau}$
is a direct summand of $[\lambda]_{\tau}\otimes[\lambda]_{\tau}$. 
Now $\chi([\mu]_{\tau})$ is given by the evaluation of the polynomial 
$p:=\frac{1}{45}(T-8)(T-6)(T-4)(T-2)(T-1)T$ at $\tau={\frac{7+i\sqrt{11}}{2}}$. 
Hence it is enough to compute the remainder of the division of $p$ by the minimum polynomial of $\tau$ over $\QQ$, 
that is $T^2-7T+15$, which is $\frac{1}{45}(135T-1080)$. 
Therefore  $\chi([\mu]_{\tau})=\frac{1}{45}(135\cdot\tau-1080)\in\CC\setminus\ZZ$.
\end{itemize}

\section*{Acknowledgements} 
The author wishes to express his gratitude to Uwe Jannsen for lots of useful conversations 
during a postdoc stay at the math department of the University of Regensburg. 
He is also indebted to Bruno Kahn, Carlo Mazza and Stefano Vigni for several suggestions and remarks.

\bibliographystyle{amsalpha}
\bibliography{Motives1}

\end{document}